\documentclass[11pt,oneside,a4paper]{article}
\usepackage{amsmath,amssymb,amsthm,epsfig,color,graphicx}
\usepackage{hyperref}
\usepackage{vmargin}
\usepackage[english,brazil]{babel}
\usepackage[latin1]{inputenc}
\usepackage{color}
\usepackage{hyperref}
\usepackage{pdfsync}
\newcommand\norm[1]{\left\lVert#1\right\rVert}
\newcommand{\R}{\mathbb{R}}

\newcommand{\normd}[2]{\left\langle  #1 \right\rangle_{#2}}
\newcommand{\noi}{\noindent}
\newenvironment{prova}{\noindent\normalsize {\sc Proof}.}{{\hfill \rule{1mm}{3mm}}}

\newtheorem{obs}{Remark}[section]
\newtheorem{lema}{Lemma}[section]

\newtheorem{teorema}{Theorem}[section]

\begin{document}

\title{Infinitely many solutions for a H\'enon-type system\\ 
in hyperbolic space}

\author{
\textbf{Patr\'icia Leal da Cunha  \footnote{\textit{E-mail addresses}:
		patcunha80@gmail.com}}\\ {\small\it Departamento de Tecnologia e Ci\^encia de Dados, FGV-SP, Brazil,} \\ 
	\\
\textbf{Fl\'{a}vio Almeida Lemos  \footnote{\textit{E-mail addresses}:
flavioal@gmail.com}}\\ 
{\small\it Departamento de Matem\'{a}tica, UFOP-MG, Brazil. }}

\date{}
\maketitle
\markboth{abstract}{abstract}
\addcontentsline{toc}{chapter}{abstract}

\noi ABSTRACT. This paper is devoted to study the semilinear elliptic system of H\'{e}non-type 
\begin{equation*}
\begin{cases}  -\Delta_{\mathbb{B}^{N}}u= K(d(x))Q_{u}(u,v), \\  -\Delta_{\mathbb{B}^{N}}v= K(d(x))Q_{v}(u,v),  \\
u, v\in H_{r}^{1}(\mathbb{B}^{N}),\, N\geq 3
\end{cases}
\end{equation*}
in the hyperbolic space $\mathbb{B}^{N}$, where $H_{r}^{1}(\mathbb{B}^{N})=\{u\in H^1(\mathbb{B}^N): u\, \text{is radial}\}$ and 
$-\Delta_{\mathbb{B}^{N}}$ denotes the Laplace-Beltrami operator on $\mathbb{B}^N$, $d(x)=d_{\mathbb{B}^{N}}(0,x)$,    
$Q \in C^{1}(\mathbb{R}\times \mathbb{R},\mathbb{R})$ is p-homogeneous 
function and $K\geq0 $ 
is a continuous function. 
We prove a compactness result and together with the Clark's theorem we establish the existence of infinitely many solutions.

\vspace{.5cm}


\section{Introduction and main result}

\noi This article concerns the existence of infinitely many solutions for the following semilinear elliptic system 
of H\'enon type in hyperbolic space
\begin{equation*}\label{problem1}
\begin{cases}  -\Delta_{\mathbb{B}^{N}}u= K(d(x))Q_{u}(u,v) \\  -\Delta_{\mathbb{B}^{N}}v= K(d(x))Q_{v}(u,v) \\  
u, v\in H_{r}^{1}(\mathbb{B}^{N}), \, N\geq 3
\tag{$\mathcal{H}$}
\end{cases}
\end{equation*}
where $\mathbb{B}^N $ is the Poincar\'e ball model for the hyperbolic space, $H_r^1(\mathbb{B}^N)$ denotes the sobolev space of radial 
$H^1(\mathbb{B}^N)$ function, $r=d(x)=d_{\mathbb{B}^{N}}(0,x)$, $\Delta_{\mathbb{B}^{N}}$ is the Laplace-Beltrami type operator on 
$\mathbb{B}^{N}$.

We assume the following hypothesis on $K$ and $Q$
\begin{itemize}
\item [$(K_1)$] $K\geq0$ is a continuous function with $K(0)=0$ and $K\neq0$  in $\mathbb{B}^{N}\backslash \{0\}$.
\item [$(K_2)$] $K=O(r^{\beta})$ as $r \to 0$ and $K=O(r^{\beta})$  as $r \to \infty$, for some $\beta >0$.
\item [$(Q_1)$] $Q \in C^{1}(\mathbb{R}\times \mathbb{R},\mathbb{R})$ is such that $Q(-s,t)= Q(s,-t)=Q(s,t)$, 
$Q(\lambda s, \lambda t)=\lambda^{p} Q(s,t)$ (Q is $p$ - homogeneous), $\forall \lambda \in \mathbb{R}$ and $p\in (2,\delta)$, 
where
\[
\delta = \begin{cases}
\dfrac{2N + 2\beta}{N-2 } & \text{if } N-2>0  \\
\infty & \text{if } \text{otherwise} \\
\end{cases}
\]
\item [$(Q_2)$] There exist $C,C_1 , C_2 >0$ such that $Q (s,t)\leq C(s^{p} + t^{p})$,  $Q_s (s,t)\leq C_1 s^{p-1}$ and 
$Q_t (s,t)\leq C_2 t^{p-1}$,\,\,  $\forall s,t \geq 0$.
\item [$(Q_3)$] There exists $C_3 >0$ such that $C_{3}(|s|^{p}+ |t|^{p})\leq Q(s,t))$ with $p\in (2,\delta)$. 
\end{itemize}

In the past few years the prototype problem
\begin{equation*}
 -\Delta_{\mathbb{B}^N}u = d(x)^{\alpha}|u|^{p-2}u, \quad u\in H_r^1(\mathbb{B}^N)
\end{equation*}
has been attracted attention.  Unlike the corresponding problem in the Euclidean space $\mathbb{R}^N$, He in \cite{HE} proved
the existence of a positive solution to the above problem over the range $p\in(2,\frac{2N+2\alpha}{N-2})$ in the hyperbolic space. 
More precisely, she explored the Strauss radial estimate for hyperbolic space together with the Mountain Pass Theorem.
In a subsequent paper \cite{HE-critical}, she proved the existence of at least one non-trivial positive solution for the critical H\'enon
equation 
\begin{equation*}
 -\Delta_{\mathbb{B}^N}u = d(x)^{\alpha}|u|^{2^*-2}u+\lambda u, \quad u\geq 0,\quad u\in H_0^1(\Omega'),
\end{equation*}
provided that $\alpha\rightarrow 0^+$ and for a suitable value of $\lambda$, where $\Omega'$ is a bounded domain in hyperbolic space 
$\mathbb{B}^N$. Finally, by working in the hole hyperbolic space $\mathbb{H}^N$, He \cite{He-sistema} considered the following 
Hardy-H\'enon type system
\begin{equation*}
\begin{cases}  
-\Delta_{\textbf{H}^{N}}u= d_b(x)^\alpha |v|^{p-1}v \\  -\Delta_{\textbf{H}^{N}}v= d_b(x)^\beta |u|^{q-1}u, \\  
\end{cases}
\end{equation*}
for $\alpha,\,\beta\in\mathbb{R}$, $N>4$ and obtained infinitely many non-trivial radial solutions.

We would like to mention the paper of Carri\~ao, Faria and Miyagaki \cite{FOC} where they extended He's result by considering
a general nonlinearity 
\begin{equation}\label{CFM}
\begin{cases}  
-\Delta_{\mathbb{B}^{N}}^\alpha u= K(d(x)) f(u) \\ u\in H_r^1(\mathbb{B}^N). \\  
\end{cases}
\end{equation}
They were able to prove the existence of at least one positive solution through a compact Sobolev embedding with the Mountain Pass
Theorem.

In this paper we investigate the existence of infinitely many solutions by considering a gradient system that generalizes
problem (\ref{CFM}). We cite \cite{BF,Chipot,DG2,FD,FF,TB,Wang-Xu-Zhang} for related gradient systems problems.
In order to obtain our result, we applied the Clark's theorem \cite{CK,DC} and get inspiration on the nonlinearities condition employed by Morais Filho and Souto \cite{FS} in a p-laplacian system defined on a bounded domain in $\mathbb{R}^N$.

As regarding the difficulties, many technical difficulties arise when working on $\mathbb{B}^N$, which is a non compact manifold. 
This means that the the embedding $H^1(\mathbb{B}^N)\hookrightarrow L^p(\mathbb{B}^N)$ is not compact for $2\leq p\leq \frac{2N}{N-2}$
and the functional related to the system \ref{problem1} cannot satisfy the $(PS)_c$ condition for all $c>0$.

We also point out that since the weight function $d(x)$ depends on the Riemannian distance $r$ from a pole o, we have some
difficulties in proving that
$$\int_{\mathbb{B}^{N}}d(x)^{\beta}\Big( |u(x)|^p + |v(x)|^p \Big) \, dV_{\mathbb{B}^{N}} <\infty, \quad \forall (u,v)\in H^1(\mathbb{B}^N)\times H^1(\mathbb{B}^N)
$$
leading to a great effort in proving that the Euler-Lagrange functional associated is well defined.

To overcome these difficulties we restrict ourselves to the radial functions.

Our result is
\begin{teorema}\label{maintheorem}
Under hypotheses $(K_1)$-$(K_2)$ and $(Q_1)$-$(Q_3)$, the problem (\ref{problem1}) has infinitely many solutions.
\end{teorema}


\section{Preliminaries}

Throughout this paper, $C$ is a positive constant which may change from line to line.

The Poincar\'e ball for the hyperbolic space is
\[
\mathbb{B}^{N}=\{x \in \mathbb{R}^{N}|\, |x|<1\}
\]
endowed with Riemannian metric $g$ given by $g_{i,j}=(p(x))^2 \delta_{i,j} $ where $p(x)= \dfrac{2}{1-|x|^2}$. 
We denote the hyperbolic volume by $dV_{\mathbb{B}^{N}}=(p(x))^N
dx$. The geodesic distance from the origin to $x \in \mathbb{B}^{N}$ is given by
\[
d(x):=d_{\mathbb{B}^{N}}(0,x)=\int_{0}^{|x|}\dfrac{2}{1-s^2}ds=\log\left(\dfrac{1+ |x|}{1-|x|}\right).
\]
The hyperbolic gradient and the Laplace-Beltrami operator are
\[
-\Delta_{\mathbb{B}^{N}}u=-(p(x))^{-N}div(p(x)^{N-2}\nabla u),  \quad
\nabla_{\mathbb{B}^{N}}u=\dfrac{\nabla u}{p(x)}
\]
\noindent where $ H^{1}(\mathbb{B}^{N})$ denotes the Sobolev space on $\mathbb{B}^{N}$ with the metric $g$. $\nabla$ and $div$ denote the 
Euclidean gradient and divergence in $\R^N$, respectively.

Let  $H^{1}_{r}(\mathbb{B}^{N})=\{ u \in H^{1}(\mathbb{B}^{N}):u \text{  is radial}  \}$.

We shall find weak-solutions of problem (\ref{problem1}) in the space
\[
H= H_{r}^{1}(\mathbb{B}^{N})\times H_{r}^{1}(\mathbb{B}^{N})
\]
endowed with the norm 
\[
\norm{(u,v)}^2= \int_{\mathbb{B}^{N}}\Big(\norm{\nabla_{\mathbb{B}^{N}}u}_{\mathbb{B}^{N}}^{2}+ 
\norm{\nabla_{\mathbb{B}^{N}}v}_{\mathbb{B}^{N}}^{2}\Big)\,dV_{\mathbb{B}^{N}}.
\]

One can observe that system (\ref{problem1}) is formally derived as the Euler-Lagrange equation for the functional
\begin{equation*}
I(u,v)= \dfrac{1}{2}\int_{\mathbb{B}^{N}}
(\norm{\nabla_{\mathbb{B}^{N}}u}_{\mathbb{B}^{N}}^{2}+ \norm{\nabla_{\mathbb{B}^{N}}v}_{\mathbb{B}^{N}}^{2})\,dV_{\mathbb{B}^{N}}-
\int_{\mathbb{B}^{N}}K(d(x))Q(u,v)\,dV_{\mathbb{B}^{N}}.
\end{equation*}

We endowed the norm for $L^{p}(\mathbb{B}^{N})\times L^{p}(\mathbb{B}^{N})$ as
\[
\norm{(u,v)}_{p}^{p}= \int_{\mathbb{B}^{N}}(|u|^{p}+ 
|v|^{p})\,dV_{\mathbb{B}^{N}}.
\]

To solve this problem we need the following lemmas.
\begin{lema}\label{lemafoc1}
The map $(u,v) \longmapsto (d(x)^m u,d(x)^m v) $ from $H= H_{r}^{1}
(\mathbb{B}^{N})\times H_{r}^{1}(\mathbb{B}^{N})$ to $L^{p}(\mathbb{B}^{N})\times L^{p}(\mathbb{B}^{N})$ 
is continuous for $ p \in (2,\tilde{m})$, where $m>0$ and
$$
\tilde{m} = \begin{cases}
\dfrac{2N}{N-2-2m  } & \text{if} \,\, m < \dfrac{N-2   }{2}, \\
\hspace{0.8cm}\infty & \text{if}\,\, \text{otherwise}. \\
\end{cases}
$$
\end{lema}

\begin{prova}
In \cite[Lemma 2.2]{HE} it has been proved that the map $u \longmapsto d(x)^m u $ from $H_{r}^{1}(\mathbb{B}^{N})$ to $L^{p}(\mathbb{B}^{N})$ 
is continuous for $ p \in (2,\tilde{m})$. Therefore
$\norm{d(x)^{m}u}_{p}\leq C \norm{u}_{H^{1}_{r}}$ and $\norm{d(x)^{m}v}_{p}\leq C \norm{v}_{H^{1}_{r}}$. 
Hence,
\begin{eqnarray*}
\Big(\norm{d(x)^{m}u}_{p}^{2}+ \norm{d(x)^{m}v}_{p}^{2}\Big)^{\frac{1}{2}}\leq C \norm{(u,v)}.
\end{eqnarray*}

Now observe that
\begin{eqnarray*}
&&\norm{d(x)^{m}u}_{p} + \norm{(d(x)^{m}v}_{p}=
\Big[\Big(\norm{d(x)^{m}u}_{p} + \norm{d(x)^{m}v}_{p}\Big)^{2}\Big]^{\frac{1}{2}}\\
&&= \Big(\norm{d(x)^{m}u}_{p}^{2} +2\norm{d(x)^{m}u}_{p} \norm{d(x)^{m}v}_{p} + \norm{d(x)^{m}v}_{p}\Big)^{\frac{1}{2}}.
\end{eqnarray*}

Applying the Cauchy's inequality $ab\leq \frac{a^2 + b^2}{2}$, we get 
\begin{eqnarray*}
\norm{d(x)^{m}u}_{p} + \norm{d(x)^{m}v}_{p} \leq \sqrt{2}
\Big(\norm{d(x)^{m}u}_{p}^{2}+ \norm{d(x)^{m}v}_{p}^{2}\Big)^{\frac{1}{2}}
\end{eqnarray*}

By the subaditivity,
\begin{eqnarray*}
\Big(\norm{d(x)^{m}u}_{p}^{p}+ \norm{d(x)^{m}v}_{p}^{p}\Big)^{\frac{1}{p}}\leq \norm{d(x)^{m}u}_{p} + \norm{d(x)^{m}v}_{p}
\end{eqnarray*}

Therefore,
\begin{eqnarray*}
\norm{(d(x)^{m}u,d(x)^{m}v)}_{p}\leq C \norm{(u,v)}
\end{eqnarray*}
and the lemma holds.
\end{prova}

\begin{obs}\label{obs1}
From the previous lemma, there exists a positive constant $C>0$, such that 
\begin{eqnarray*}
&&\int_{\mathbb{B}^{N}}d(x)^{\beta}|u(x)|^p dV_{\mathbb{B}^{N}} + \int_{\mathbb{B}^{N}}d(x)^{\beta}|v(x)|^p dV_{\mathbb{B}^{N}}\\
&& \leq C \left( \int_{\mathbb{B}^{N}}\norm{\nabla_{\mathbb{B}^{N}}u}_{\mathbb{B}^{N}}^{2}dV_{\mathbb{B}^{N}} + 
\int_{\mathbb{B}^{N}}\norm{\nabla_{\mathbb{B}^{N}}v}_{\mathbb{B}^{N}}^{2}dV_{\mathbb{B}^{N}} \right)^{\frac{p}{2}},
\end{eqnarray*}
where $m=\frac{\beta}{p}$ and $2<p<\frac{2N}{N-2-2(\frac{\beta}{p})}$, that is, $2<p<\delta$.
\end{obs}

\begin{lema}\label{teoremafoc1}
The map $(u,v) \longmapsto (d(x)^m u,d(x)^m v) $ from $H= H_{r}^{1}(\mathbb{B}^{N})\times H_{r}^{1}(\mathbb{B}^{N})$ to 
$L^{p}(\mathbb{B}^{N})\times L^{p}(\mathbb{B}^{N})$ is compact for $ p \in (2,\tilde{m})$, where $m>0$ and		
\[
\tilde{m} = \begin{cases}
\dfrac{2N}{N-2-2m  } & \text{if } m< \dfrac{N-2   }{2} \\
\hspace{0.8cm}\infty & \text{if } \text{otherwise} \\
\end{cases}
\]
\end{lema}

\begin{prova}
Let $(u_{n}, v_{n}) \in H $ be a bounded sequence. Then up to a subsequence, if necessary, we may assume that
\[
(u_{n}, v_{n})\rightharpoonup (u,v)
\]
	
It is easy to sea that $u_n \rightharpoonup u$ and $v_n \rightharpoonup v$ in $ H_{r}^{1}(\mathbb{B}^{N})$. 

We will use the same calculus used by Haiyang He  \cite{HE} (page $26$). We want  to show that  
\[
\lim_{n\to\infty}\int_{\mathbb{B}^{N}}d(x)^{mq}|u_{n}(x)|^p dV_{\mathbb{B}^{N}} = \int_{\mathbb{B}^{N}}d(x)^{mp}|u(x)|^p dV_{\mathbb{B}^{N}}, 
\]
\[
\lim_{n\to\infty}\int_{\mathbb{B}^{N}}d(x)^{mq}|v_{n}(x)|^p dV_{\mathbb{B}^{N}} = \int_{\mathbb{B}^{N}}d(x)^{mp}|v(x)|^p dV_{\mathbb{B}^{N}}.
\]
Let 	$u\in H_{r}^{1}(\mathbb{B}^{N})$, then by Haiyang He \cite{HE} we have
\begin{eqnarray*}
|u_{n}(x)|&&\leq \dfrac{1}{\sqrt{\omega_{n-1}(N-2)}}\left(\frac{1- |x|^2}{2}
\right)^{\frac{N-2}{2}}\dfrac{1}{|x|^{\frac{N-2}{2}}}\norm{u_{n}}_{H_{r}^{1}(\mathbb{B}^{N})},\\
|u_{n}(x)|&&\leq \dfrac{1}{\sqrt{\omega_{n-1}}}\left(\frac{1- |x|^2}{2}\right)^{\frac{N-1}{2}}\dfrac{1}{|x|^{\frac{N}{2}}}
\norm{u_{n}}_{H_{r}^{1}(\mathbb{B}^{N})}.
\end{eqnarray*}
Since $\{|x|\leq \frac{1}{2}\}$,  $\ln\dfrac{1+ |x|}{1-|x|}\leq \dfrac{2r}{1-r^2}$ and $2<p< \tilde{m}$, we have

\begin{eqnarray*}
d(x)^{mp}|u|^{p}&&\leq C \left(\ln\dfrac{1+ |x|}{1-|x|}\right)^{mp}\left(\frac{1- |x|^2}{2}\right)^{p\frac{N-2}{2}}
\left(\dfrac{1}{|x|^{\frac{N-2}{2}}}\right)^{p}\\
	&& \leq C\left(\dfrac{2|x|}{1-|x|^2}\right)^{mp}\left(\frac{1- |x|^2}{2}\right)^{p\frac{N-2}{2}}\left(\dfrac{1}
	{|x|^{\frac{N-2}{2}}}\right)^{p} \equiv h_{1}.
\end{eqnarray*}

Set 
\[
g_{1}(x) = \begin{cases}
h_{1}(x) & \text{if } 0\leq |x| < \frac{1}{2} \\
0 & \text{if }\frac{1}{2}\leq |x| <1,  \\
\end{cases}
\]
then
\begin{eqnarray*}
\int_{\mathbb{B}^{N}}g_1 \,dV_{\mathbb{B}^{N}}&&=\int_{0}^{\frac{1}{2}}\left(\dfrac{2r}{1-r^2}\right)^{mp}
\left(\frac{1- r^2}{2}\right)^{p\frac{N-2}{2}}\left(\dfrac{1}{r^{\frac{N-2}{2}}}\right)^{p}r^{N-1}\left(\frac{2}{1- r^2}\right)^{N}dr \\
&&\leq C\int_{0}^{\frac{1}{2}}\left(\dfrac{2r}{1-r^2}\right)^{mp}\left(\frac{1- r^2}{2}\right)^{(p\frac{N-2}{2} - N)}r^{N-1 - p\frac{N-2}{2}}dr\\
&&\leq C \int_{0}^{\frac{1}{2}}r^{mp +N -1 -p\frac{N-2}{2}}\,dr < \infty.
\end{eqnarray*}
 
Since $\{|x|> \frac{1}{2}\}$ and $2<p< \tilde{m}$, we have
\begin{eqnarray*}
	d(x)^{mp}|u|^{p}&&\leq C \left(\ln\dfrac{1+ |x|}{1-|x|}\right)^{mp}\left(\frac{1- 
	|x|^2}{2}\right)^{p\frac{N-1}{2}}\left(\dfrac{1}{|x|^{\frac{N}{2}}}\right)^{p}\\
	&& \leq C\left(\ln\dfrac{1+ |x|}{1-|x|}\right)^{mp}\left(\frac{1- |x|^2}{2}
	\right)^{p\frac{N-1}{2}}\left(\dfrac{1}{|x|^{\frac{N}{2}}}\right)^{p} \equiv h_{2}.
\end{eqnarray*}
Set 
 \[
g_{2}(x) = \begin{cases}
0 & \text{if } 0\leq |x| < \frac{1}{2} \\
h_{2}(x) & \text{if }\frac{1}{2}\leq |x| <1,  \\
\end{cases}
\]
then
\begin{eqnarray*}
	\int_{\mathbb{B}^{N}}g_2 \,dV_{\mathbb{B}^{N}}&&=\int_{\frac{1}{2}}^{1}\left(\ln\dfrac{1+ r}{1-r}\right)^{mp}
	\left(\frac{1- r^2}{2}\right)^{p\frac{N-1}{2}}\left(\dfrac{1}{r^{\frac{N}{2}}}\right)^{p}r^{N-1}\left(\frac{2}{1- r^2}\right)^{N}dr \\
	&&\leq C\int_{\frac{1}{2}}^{1}\left(\ln\dfrac{1+ r}{1-r}\right)^{mp}\left(\frac{1- r^2}{2}\right)^{(p\frac{N-1}{2} - N)}r^{N-1 - \frac{N}{2}p}dr\\
	&&\leq \int_{\ln 3}^{\infty}s^{mp}\left(\dfrac{2e^{s}}{(e^{s} +1)^2}\right)^{\frac{N-1}{2}p-N +1}ds< \infty.
\end{eqnarray*}
Hence, we have 
\[
|d(x)^{mq}u_{n}(x)^{q}|\leq g_{1}(x) + g_{2}(x).
\]
By the Dominated convergence theorem, we obtain
\[
\lim_{n\to\infty}\int_{\mathbb{B}^{N}}d(x)^{mq}u_{n}^{p}(x) dV_{\mathbb{B}^{N}} = \int_{\mathbb{B}^{N}}d(x)^{mp}u^{p}(x) dV_{\mathbb{B}^{N}}. 
\]
In the same way we conclude that
\[
\lim_{n\to\infty}\int_{\mathbb{B}^{N}}d(x)^{mq}v_{n}^{p}(x) dV_{\mathbb{B}^{N}} = \int_{\mathbb{B}^{N}}d(x)^{mp}v^{p}(x) dV_{\mathbb{B}^{N}} 
\]
and the Lemma holds.
\end{prova}


\section{Proof of Theorem \ref{maintheorem}}

The Clark's theorem is one of the most important results in critical point theory (see \cite{CK}). 
It was successfully applied to sublinear elliptic problems with symmetry and the existence of infinitely many solutions around the $0$ was shown.
 
In order to state the Clark's theorem, we need some terminologies. 

Let $(X,\norm{\cdot}_{X})$ be a Banach Space and $\mathcal{I} \in C^{1}(X,\R)$.
\begin{itemize}
\item[(i)] For $c \in \R$ we say that $\mathcal{I}(u)$ satisfies $(PS)_{c}$ condition if any sequence 
$(u_j)_{j=1}^{\infty} \subset X$ such that $\mathcal{I}(u_{j})\to c$ and $\norm{\mathcal{I}'(u_j)}\to 0$ has a convergent subsequence.
\item[(ii)] Let  $S$ be a symmetric and closed set family in $X\diagdown \{0\}$. 
For  $A \in S$, the genus  $\gamma(A)=\min\{n \in \mathbb{N}:\phi \in C(A, \mathbb{R}^n \{0\}) \text{ is odd}\}$. 
If there is no such natural, we set $\gamma(A)= \infty$.
\item[(iii)] Let $\Omega$ be  a open and bounded set, $0\in\Omega$  in $\R^n$. If $A\in S$ is such that there exists a odd homeomorphism 
function from $A$ to $\partial \Omega$, then $\gamma(A)=n$.
\end{itemize}
 
\begin{teorema}[Clark's Theorem]\label{index}
Let $\mathcal{I} \in C(X, \R)$ be an even function, bounded from below, with $\mathcal{I}(0)=0$ and there exists a compact, symmetric set $K \in S$ 
such that $\gamma(K)=k$ and $\sup_{K}\mathcal{I}<0$. Then $I$ has least $k$ distinct pairs of critical points.
\end{teorema}

The proof of Theorem \ref{maintheorem} is made by using Theorem \ref{index}.

The (\ref{problem1}) system are the Euler-Lagrange equations related to the functional
\begin{equation}\label{I2}
I(u,v)= \dfrac{1}{2}\int_{\mathbb{B}^{N}}(\norm{\nabla_{\mathbb{B}^{N}}u}_{\mathbb{B}^{N}}^{2}+ 
\norm{\nabla_{\mathbb{B}^{N}}v}_{\mathbb{B}^{N}}^{2})dV_{\mathbb{B}^{N}}-
\int_{\mathbb{B}^{N}}K(d(x))Q(u,v)dV_{\mathbb{B}^{N}}
\end{equation}
which is $C^1$ on $H$. 

The functional $I$ is not bounded from below,
therefore, we can't apply the  Clark's technique  for this functional. 

In order to overcome this difficulty we consider the auxiliary functional
\begin{equation} \label{mf1}
J(u,v)= \left(\int_{\mathbb{B}^{N}}(\norm{\nabla_{\mathbb{B}^{N}}u}_{\mathbb{B}^{N}}^{2}+
 \norm{\nabla_{\mathbb{B}^{N}}v}_{\mathbb{B}^{N}}^{2})\, dV_{\mathbb{B}^{N}} \right)^{p-1} -
 \int_{\mathbb{B}^{N}}K(d(x))Q(u,v)\,dV_{\mathbb{B}^{N}},
\end{equation}
where $p \in(2, \delta)$, while for $J'$ we have, $\forall (\phi,\psi)\in H$
\begin{eqnarray} \label{mf2}
J'(u,v)(\phi, \psi) &=& (2p-2)\norm{(u,v)}^{2p-4}\int_{\mathbb{B}^{N}}
(\normd{\nabla_{\mathbb{B}^{N}}u,\nabla_{\mathbb{B}^{N}}\phi}{\mathbb{B}^{N}}+
\normd{\nabla_{\mathbb{B}^{N}}v,\nabla_{\mathbb{B}^{N}}\psi}{\mathbb{B}^{N}})\,dV_{\mathbb{B}^{N}}\nonumber \\
&&- \int_{\mathbb{B}^{N}}K(d(x)) (\phi Q_{u}(u,v) + \psi Q_{v}(u,v)) \,dV_{\mathbb{B}^N}
\end{eqnarray}

We will show that the  set of critical points of $J$ is related to a set of critical points of $I$ and $J$ 
satisfies the conditions of Theorem \ref{index}.

The proof of Theorem \ref{maintheorem} is divided into several lemmas.

\begin{lema}\label{sol1}
If $(u,v)\in H$, $(u,v) \neq (0,0)$ is a critical point for $J$, then 
\[(w,z)=\left(\dfrac{u}{[(2p-2)\norm{(u,v)}^{2p-4}]^{\frac{1}{p-2}}}, \dfrac{v}{[(2p-2)\norm{(u,v)}^{2p-4}]^{\frac{1}{p-2}}}\right)\]
is a critical point for $I$.
\end{lema}

\begin{prova}
Note that $(u,v)\neq (0,0)$ is a critical point for $J$ if, and only if, $(u,v)$ is a weak solution to the problem 
\begin{equation*}\label{problem2}
\begin{cases}  -(2p-2)\norm{(u,v)}^{2p-4}\Delta_{\mathbb{B}^{N}}u= K(d(x))Q_{u}(u,v) \\ 
 -(2p-2)\norm{(u,v)}^{2p-4}\Delta_{\mathbb{B}^{N}}v= K(d(x))Q_{v}(u,v) \\  u, v\in H_{r}^{1}(\mathbb{B}^{N}), \, N\geq 3
	\tag{$\mathcal{S}$}
\end{cases}
\end{equation*}

Define $\displaystyle\lambda (\norm{(u,v)})=\Big[(2p-2)\norm{(u,v)}^{2p-4}\Big]^{\frac{-1}{p-2}}$, 
then $(w,z)=\lambda (\norm{(u,v)}) (u,v) $. 

Using the $p-1$-homogeneity condition of $Q_{u}(u,v)$ and $Q_{u}(u,v)$, observe that
\begin{eqnarray*}
&& - \Delta_{\mathbb{B}^N}w-K(d(x))Q_{u}(w,z)\\
&& = -\lambda (\norm{(u,v)})\Delta_{\mathbb{B}^N}u - (\lambda (\norm{(u,v)}))^{p-1}K(d(x))Q_{u}(u,v)\\
&& = -\lambda (\norm{(u,v)}) K(d(x)) Q_{u}(u,v)\left( 
(\lambda (\norm{(u,v)}))^{p-2} - \dfrac{1}{(2p-2)\norm{(u,v)}^{2p-4}} \right)
\end{eqnarray*}
and
\begin{eqnarray*}
	&& - \Delta_{\mathbb{B}^N}z-K(d(x))Q_{v}(w,z)\\
	&& = -\lambda (\norm{(u,v)})\Delta_{\mathbb{B}^N}v - (\lambda (\norm{(u,v)}))^{p-1}K(d(x))Q_{v}(u,v)\\
	&& = -\lambda (\norm{(u,v)}) K(d(x)) Q_{v}(u,v)\left( 
	(\lambda (\norm{(u,v)}))^{p-2} - \dfrac{1}{(2p-2)\norm{(u,v)}^{2p-4}} \right)
\end{eqnarray*}
Hence $(w,z)$ is a weak solution for problem (\ref{problem1}) and so, a critical point for $I$.
\end{prova}

\begin{lema}\label{sol2}
$J(u,v)$ is bounded from below and  satisfies the $(PS)_c$ condition.
\end{lema}

\begin{prova}
From  $(K_1)$-$(K_2)$, $(Q_2)$-$(Q_3)$ and Remark \ref{obs1}
\begin{eqnarray*}
J(u,v)&=& 
 \left(\int_{\mathbb{B}^{N}}\norm{\nabla_{\mathbb{B}^{N}}u}_{\mathbb{B}^{N}}^{2}+
  \norm{\nabla_{\mathbb{B}^{N}}v}_{\mathbb{B}^{N}}^{2}\, dV_{\mathbb{B}^{N}} \right)^{p-1} -
\int_{\mathbb{B}^{N}}K(d(x))Q(u,v)\,dV_{\mathbb{B}^{N}}\\
&\geq &  \norm{(u,v)}^{2p-2} - C\int_{\mathbb{B}^{N}}d(x)^{\beta}(|u|^{p} +|v|^{p})\,dV_{\mathbb{B}^{N}}\\
&\geq & \norm{(u,v)}^{2p-2} - \norm{(u,v)}^{p}, 
\end{eqnarray*}
\noindent so that $J(u,v)$ is bounded from below.

Let $(u_n.v_n) \in H$ be such that $|J(u_n,v_n)|\leq C$ with $C\in\R^+$, $J'(u_n,v_n)\to 0$. Since
\[
C\geq J(u_n ,v_n)\geq \norm{(u_n,v_n)}^{2p-2} - \norm{(u_n,v_n)}^{p},
\]
we conclude that $\norm{(u_n,v_n)} $  is bounded. So, there exists $(u,v) \in H$ such that, 
passing to a subsequence if necessary,
$$
(u_n, v_n) \rightharpoonup (u,v), \quad \mbox{as} \,\, n \rightarrow \infty.
$$

From the Embedding Lemma \ref{teoremafoc1},
\[\int_{\mathbb{B}^N}(d(x))^{\beta}(|u_n|^p + |v_n|^p) \,dV_{\mathbb{B}^N}\longrightarrow 
\int_{\mathbb{B}^N}(d(x))^{\beta}(|u|^p + |v|^p )\, dV_{\mathbb{B}^N}. 
\]
and by $(K_2)-(Q_2)$, we infer that
$$
|K(d(x)(u_nQ_{u}(u_n,v_n) +v_nQ_{v}(u_n,v_n) )|\leq C(d(x))^{\beta}(|u_n |^{p}+ v_n |^{p}).
$$

Therefore, by the  Lebesgue Dominated Convergence Theorem,  
\[
\int_{\mathbb{B}^{N}}K(d(x))( Q_{u}(u_n , v_ n)u_n +Q_{v}(u_n,v_n)v_n )\, dV_{\mathbb{B}^{N}}\rightarrow 
\int_{\mathbb{B}^{N}}K(d(x))(Q_{u}(u,v)u +Q_{v}(u,v )v)\, dV_{\mathbb{B}^{N}}.
\]

Since $J'(u,v)(u,v)=0$ and $J'(u_n,v_n)(u_n, v_n)=o_{n}(1)$ as $n\to \infty$, we have
\begin{eqnarray*}
&&(2p-2)\Big (  \norm{(u_n,v_n)}^{2p-2}-\norm{(u,v)}^{2p-2}\Big )=J'(u_n,v_n)(u_n, v_n)-J'(u,v)(u,v)+\\
&&+\int_{\mathbb{B}^{N}}K(d(x))\Big(Q_{u}(u_n , v_ n)u_n +Q_{v}(u_n,v_n)v_n \Big)\,dV_{\mathbb{B}^N}+\\
&&-\int_{\mathbb{B}^{N}}K(d(x))\Big(Q_{u}(u,v)u +Q_{v}(u,v )v\Big)\,dV_{\mathbb{B}^N}= o_n(1)
\end{eqnarray*}
\noindent then, 
$\norm{(u_n,v_n)}\rightarrow \norm{(u,v)}$. Therefore, 
$$
(u_n,v_n)\rightarrow (u,v), \quad \mbox{as} \,\, n \rightarrow \infty, \,\, \mbox{in} \,\, H
$$
\end{prova}

The next lemma ends the proof of the Theorem \ref{maintheorem}.

\begin{lema}\label{sol3}
	Given $k \in \mathbb{N}$, there exists a compact and symmetric set $K \in H$ such that $\gamma(K)=k$ and $\sup_{K}J <0$.
\end{lema}

\begin{prova}
Let $X_k \subset H$ be a subspace of dimension $k$. Consider the following norm in $X_k$
	\[
	\norm{(u,v)}_{X_k}=\left( \int_{\mathbb{B}^N}K(d(x))(|u|^p + |v|^p) dV_{\mathbb{B}^N}\right)^{\frac{1}{p}}.
	\]
	
	 Since $X_k \subset H$ has finite dimension, there exists $a>0$ such that
	\[
	a\norm{(u,v)}_{X_k}\leq \norm{(u,v)} \leq \frac{1}{a}\norm{(u,v)}_{X_k}, \quad \forall (u,v)\in X_k.
	\]
	
Therefore, we obtain from $(Q_3)$ that
	\[
	J(u,v)\leq \norm{(u,v)}^{2p-2}-C\int_{\mathbb{B}^{N}}K(d(x))(|u|^{p}+ |v|^{p})dV_{\mathbb{B}^N} = \norm{(u,v)}^{2p-2} - C \norm{(u,v)}_{X_k}^{p},
	\]
\noindent where $C\in\mathbb{R}$ is a positive constant. We then conclude that 
	\[
	J(u,v)\leq \norm{(u,v)}_{X_k}^{p}\left( \dfrac{\norm{(u,v)}_{X_k}^{p-2}}{a^{2p-2}}- C\right).
	\]
	
Let $A=a^{\frac{2p-2}{p-2}}$ and consider the set $K=\{ (u,v)\in X_k: \norm{(u,v)}_{X_k}= \frac{A}{2}C^{\frac{1}{p-2}} \}$, then
	\[
	J(u,v)\leq C \frac{A^p}{2^p}\left( \dfrac{1}{2^{q-2}}-1\right)<0,\quad \forall (u,v) \in K.
	\]
	
We get that $\sup_{K}J <0$, where $K\subset H$ is a compact and symmetric set such that $\gamma(K) =k$.
	
\end{prova}

Finally, from Lemmas \ref{sol2} and \ref{sol3}, Theorem \ref{index}  implies the existence of at least $k$ distinct pairs 
of critical points for the functional $J$. Since $k$ is arbitrary, we obtain infinitely many critical points for $J$.
	
In view of Lemma \ref{sol1},  we conclude that the functional $J$ possesses,
together with $I$, infinitely many critical points. 

Finally, we point out that since $H$ is a closed subspace of the Hilbert space $H^1(\mathbb{B}^N)\times H^1(\mathbb{B}^N)$, 
following some ideas in \cite{BCS,FOC}, we can conclude that $(u,v)$ is a critical point in 
$H^{1}(\mathbb{B}^{N})\times H^{1}(\mathbb{B}^{N})$.

\section{Further result}

We can apply the same method that proved Theorem \ref{maintheorem} to establish the existence of infinitely many solutions for
the following semilinear elliptic equation
\begin{equation*}\label{problem2}
\begin{cases}  -\Delta_{\mathbb{B}^{N}}^{\alpha}u= K(d(x))|u|^{p-2}u, \\  
u\in \mathbb{E}\subset H_{r}^{1}(\mathbb{B}^{N}), \, N\geq 3
\end{cases}\tag{$\mathcal{H*}$}
\end{equation*}
where $K$ satisfies $(K_1)-(K_2)$, $-\Delta_{\mathbb{B}^{N}}^{\alpha}$ is the Laplace-Beltrami type operator
\[
-\Delta_{\mathbb{B}^{N}}^{\alpha}u=-(p(x))^{-N}div(p(x)^{N-2}(d(x))^{\alpha}\nabla u)
\]
and
\begin{eqnarray*}
\mathbb{E}= \left\{ u \in H^{1}_{r}(\mathbb{B}^{N}): \norm{u}_{\mathbb{E}} =
\left( \int_{\mathbb{B}^{N}}(d(x))^{\alpha}\norm{\nabla_{\mathbb{B}^{N}}u}^2_{\mathbb{B}^{N}}dV_{\mathbb{B}^{N}}\right)^\frac{1}{2} \right\}.
\end{eqnarray*}
We obtain the following result
\begin{teorema}\label{maintheorem2}
	Under hypotheses $(K_1)$-$(K_2)$, (\ref{problem2}) equation has infinitely many solutions.
\end{teorema}

The energy functional corresponding to (\ref{problem2}) is
\begin{equation}\label{I1}
I(u)= \dfrac{1}{2}\int_{\mathbb{B}^{N}}(d(x))^{\alpha}\norm{\nabla_{\mathbb{B}^{N}}u}_{\mathbb{B}^{N}}^{2}dV_{\mathbb{B}^{N}}-
\frac{1}{q}\int_{\mathbb{B}^{N}}K(d(x))|u|^q\,dV_{\mathbb{B}^{N}}
\end{equation}
defined on $H_{r}^{1}(\mathbb{B}^{N})$

Problem (\ref{problem2}) is close related to one studied by Carri\~ao, Faria and Miyagaki \cite{FOC}. 
In \cite{FOC}, they proved that the map $u \longmapsto d(x)^m u $ from $\mathbb{E}$ to $L^{q}(\mathbb{B}^{N})$ is 
compact for $ q \in (2,\tilde{m})$, where
	\[
	\tilde{m} = \begin{cases}
	\dfrac{2N}{N-2-2m + \alpha} & \text{if } m< \dfrac{N-2 + \alpha }{2} \\
	\infty & \text{if } \text{otherwise} \\
	\end{cases}
	\]
and then there exists a positive constant $C>0$, such that
\begin{eqnarray}\label{remark}
	\int_{\mathbb{B}^{N}}d(x)^{\beta}|u(x)|^q dV_{\mathbb{B}^{N}}
	\leq C \left( \int_{\mathbb{B}^{N}}(d(x))^{\alpha}\norm{\nabla_{\mathbb{B}^{N}}u}_{\mathbb{B}^{N}}dV_{\mathbb{B}^{N}} \right)^{\frac{q}{2}},
\end{eqnarray}
by taking $m=\frac{\beta}{q}$ with $2<q<\frac{2N}{N-2-2\frac{\beta}{q}+ \alpha}$.

Using $(K_1)$-$(K_2)$ together with inequality (\ref{remark}), we get that the functional $I$ is well defined.
This functional is not bounded from below, hence, we can't apply Clark's technique \cite{CK}. 

In order to overcome this difficulty we consider the auxiliary functional
\begin{equation} \label{mf1}
\psi(u)= \left(\int_{\mathbb{B}^{N}}(d(x))^{\alpha}\norm{\nabla_{\mathbb{B}^{N}}u}_{\mathbb{B}^{N}}^{2}dV_{\mathbb{B}^{N}} \right)^{p-1} -
\int_{\mathbb{B}^{N}}K(d(x))|u|^{p}dV_{\mathbb{B}^{N}},
\end{equation}
where $p \in(2, 2_{\alpha}^{\beta})$, $u\in \mathbb{E}$ and

\begin{eqnarray} \label{mf2}
\psi'(u)v &=& (2p-2)\norm{u}_{\mathbb{E}}^{2p-4}\int_{\mathbb{B}^{N}}(d(x))^{\alpha}\normd{\nabla_{\mathbb{B}^{N}}u,
	\nabla_{\mathbb{B}^{N}}v}{\mathbb{B}^{N}}\,dV_{\mathbb{B}^N}+\nonumber \\
&&- \int_{\mathbb{B}^{N}}K(d(x))|u|^{p-2}uv \,dV_{\mathbb{B}^N}
\end{eqnarray}

We have the corresponding results of Lemmas \ref{sol1}, \ref{sol2} and \ref{sol3} for problem (\ref{problem2}). 
The  set of critical points of $\psi$ is related to a set of critical points of $I$ and $\psi$ 
satisfies the conditions of Theorem \ref{index}. 
\begin{lema}\label{soll1}
	If $u\in \mathbb{E}$, $u \neq 0$ is a critical point for $\psi$, then $v=\dfrac{u}{[(2p-2)\norm{u}_{\mathbb{E}}^{2p-4}]^{\frac{1}{p-2}}}$ 
	is a critical point for $I$,
\end{lema}
 
\begin{lema}\label{soll2}
	$\psi(u)$ is bounded from below and  satisfies the Palais-Smale condition (PS).
\end{lema}

\begin{lema}\label{soll3}
	Given $k \in \mathbb{N}$, there exists a compact and symmetric set $K \in \mathbb{E}$ such that $\gamma(K)=k$ and $\sup_{K}\psi <0$.
\end{lema}

From Lemmas  \ref{soll2} and  \ref{soll3}  and Theorem \ref{index}, we conclude that the functional $I$ possesses infinitely many critical points.

\vspace{.5cm}
{\bf Acknowledgments.} This paper was carried out while the first author was visiting the Mathematics Department of UFOP and
she would like to thank the members of DEMAT-UFOP for their hospitality. 
The authors are greatful to professor Giovany Figueiredo for helpful comments.


\end{document}